\documentclass[12pt]{article}

\input diagxy.tex
\xyoption{curve}

\newtheorem{thm}{Theorem }[section]
\newtheorem{df}[thm]{Definition }
\newtheorem{lm}[thm]{Lemma }
\newtheorem{ex}[thm]{Example}
\newtheorem{cor}[thm]{Corollary }

\newcommand{\pf}{{\noindent\bf Proof: }}

\newcommand{\qed}{\hfill\rule[-1pt]{4pt}{8pt}\par\vspace{\baselineskip}}
\newcommand{\edf}{\hfill$\diamondsuit$\end{df}}
\newcommand{\eex}{\hfill\framebox[1.2\width]{\ }\par\end{ex}}

\newcommand{\ds}{\displaystyle}

\newcommand{\vsc}[1]{\vspace{#1cm}}

\newcommand{\lra}{\longrightarrow}
\newcommand{\ra}{\rightarrow}

\newcommand{\Ra}{\Rightarrow}

\renewcommand{\t}[1]{\hbox{{\tiny $#1$}}}
\newcommand{\arr}[4]{#1\!\stackrel{\t{#2}}{\lra}_{#4}\!#3}

\newcommand{\twoarr}[5]{#1\!\stackrel{\t{#2}}{\ra}\!#3\!\stackrel{\t{#4}}{\ra}\!#5}

\newcommand{\tran}[3]{#1\stackrel{#2}{\Longrightarrow}#3}
\newcommand{\twotran}[5]{#1\stackrel{#2}{\Longrightarrow}{#3}\stackrel{#4}{\Longrightarrow}{#5}}

\newcommand{\Tran}[3]{#1\!:\!#2\!\Rightarrow\!#3}

\newcommand{\sub}[1]{_{_{#1}}}

\newcommand{\alphasingle}[2]{\alpha^{#1}_{#2}}

\newcommand{\h}{{\cal H}}
\newcommand{\q}{{\cal Q}}

\newcommand{\C}{{\cal C}}
\newcommand{\D}{{\cal D}}

\renewcommand{\L}{{\cal L}}

\newcommand{\V}{{\cal V}}
\newcommand{\X}{{\cal X}}

\newcommand{\ord}{{\cal O}}

\newcommand{\cat}{{\bf Cat}}
\newcommand{\scat}{{\bf Scat}}

\newcommand{\Ord}{{\bf Ord}}
\renewcommand{\inf}{{\bf Inf}}
\newcommand{\set}{{\bf Set}}
\newcommand{\rel}{{\bf Rel}}
\newcommand{\lax}{{\bf Lax}}
\newcommand{\pre}{{\bf Pre}}
\newcommand{\Shv}{{\bf Shv}}

\newcommand{\kar}{{\bf Kar}}
\newcommand{\skar}{{\bf Kar^*}}
\newcommand{\map}{{\bf Map}}

\newcommand{\Preshv}{\set^{\h^{op}}}

\newcommand{\shv}{{\bf SHV}(\h)}
\newcommand{\shvd}{{\bf SHV}\Big(\D(\h)\Big)}

\newcommand{\Relpreshv}{\rel^{\h^{co}}_{\Ord}}

\newcommand{\Relshv}{\rel^{\h^{co}}_{\inf}}
\newcommand{\Relshvd}{\rel^{\D(\h)^{co}}_{\inf}}

\newcommand{\meet}{\wedge}

\title{Relational-Sheaves for a Heyting Algebra}
\author{W. Dale Garraway}
\date{}

\begin{document}

\maketitle

\vspace{-1.5cm}

\begin{abstract}

We show that for a Heyting algebra $\h$, a {\it relational-presheaf} is an idempotent symmetric order-preserving lax-semifunctor.  A relational-presheaf is a {\it relational-sheaf} if it is an idempotent infima-preserving lax semifunctor.  The associated relational-sheaf functor factors through the category of sheaves for $\h$.  Using this and the appropriate comparison theorems we obtain the main result that the associated categories of relational-presheaves and relational-sheaves are each respectively equivalent to the categories of presheaves and sheaves for $\h$.

\end{abstract}

\section{Introduction}

For $\h$ a complete Heyting algebra the category of presheaves consists of contravariant functors $F:\h^{co}\ra\set$ and the transformations between them and  a sheaf is a presheaf that satisfies the appropriate gluing condition.  The relationships that exist between the category of presheaves and the category of sheaves for a complete Heyting algebra is well understood and the work here focuses on two of these relationships: The comparison theorem which is an equivalence between the category of presheaves for $\h$ and the category of sheaves on the Heyting algebra of down-closed subsets and the associated sheaf functor which converts a presheaf into a sheaf.

$$\bfig
\morphism(0,0)|a|/{@{>}@/^1.5em/}/<1000,0>[\shvd`\Preshv;\Lambda]
\morphism(0,0)|b|/{@{<-}@/_1.5em/}/<1000,0>[\shvd`\Preshv;\Gamma]
\morphism(1000,0)|a|/{@{>}@/^1.5em/}/<1000,0>[\Preshv`\shv;a_{\Shv}]
\morphism(1000,0)|b|/{@{<-}@/_1.5em/}/<1000,0>[\Preshv`\shv;]
\put(45,-3){{\tiny comparison}}
\put(59,3){{\tiny $\sim$}}
\put(173,-2){$\bot$}
\efig$$

Rosenthal\cite{ro} defined a relational-presheaf on a Heyting algebra\footnote{In fact he defined them for a supremum-enriched category (quantaloid) of which a Heyting algebra is a particular example.} $\h$, to be a lax-functor $F:\h^{co}\ra\rel$ (codomain the category of sets and relations), and a morphism of relational-presheaves is a lax-natural transformation in which each morphism is a function.  In his work the Heyting algebra is interpreted as a one object supremum enriched category as opposed to being a partial order with extra structure.  For him a relational-presheaf is said to be continuous if it preserves infima.  In Garraway\cite{wdg2}  a relational-presheaf is generalized to be a lax-semifunctor with added structure and the category of relational-sheaves ($\Relshv$) then arises naturally using the Karoubian-envelope construction and it is shown that the category of relational-sheaves is isomorphic to the category of $\q$-valued sets for a quantaloid $\q$.   In this paper we will explicitly define relational-presheaves to be order preserving lax-semifunctors with added structure ($\Relpreshv$ denotes the associated category).   With this in mind  we will construct a comparison theorem and an associated relational-sheaf functor for relational-presheaves  and relational-sheaves.   As a consequence of our constructions we will show that the categories of sheaves and presheaves are  equivalent to the categories of relational-sheaves and relational-presheaves respectively.


We will begin with a review of the basic definitions and constructions that are needed to define a complete Heyting algebra ($\h$) as a partial order with added structure.  This is the interpretation used when constructing the category of presheaves and sheaves for $\h$.  The first section finishes with a proof of the comparison theorem then relates presheaves on $\h$ with sheaves on the down-closed subsets of $\h$.    

Following this we then give a review of order-enriched semicategories and sup\-re\-mum-enriched semicategories.  We will use the fact that both a Heyting algebra and the category of sets and relations (\rel)  are supremum-enriched and are involutive when we define relational-sheaves and relational-presheaves.   Before we can get to the complete definition though we need review the concept of  the idempotent splitting completion (the Karoubian-envelope) of a category.  This will be the main tool that we will use to construct the categories of relational-presheaves and relational-sheaves 

In the next section we explicitly define the categories of relational-presheaves and relational-sheaves in two ways.  The first method will start with families of arrows as the objects and the morphisms between them and then constructs the categories of relational-presheaves and sheaves using the Karoubian-envelope. Taking this route gives us two things; it represents the construction as the coproduct completion as an enriched semicategory followed by the idempotent splitting completion, and secondly it gives us a slightly easier proof of the comparison theorem for relational-presheaves and sheaves.  Once this is done we then show that this construction is equivalent to defining the category of relational-presheaves as symmetric idempotent order-preserving lax-semifunctors and the category or relational-sheaves as symmetric idempotent infimum-preserving lax-semifunctors

We finish by  creating the main adjunction that we use to relate category of sheaves and the category of relational-sheaves.  
$$\Delta_{\inf}\dashv\Theta_{\inf}:\Relshv\ra\Preshv$$
When this is restricted to the category of sheaves for $\h$ the adjunction becomes an equivalence.  Once we have created this adjunction all of the equivalences and associated sheaf and relational-sheaf functors  fall out naturally.  These are represented in the diagram below

\vsc{-1}

$$\bfig
\square(0,500)|alla|/->`>->`>->`->/<950,500>[\shv`\Relshv`\Preshv`\Relpreshv;\Delta_{\Shv}`\Theta\Delta`\Delta\Theta`\Delta_{\pre}]
\square|allb|/`>`>`->/<950,500>[\Preshv`\Relpreshv`\shvd`\Relshvd;`\Gamma`\Psi`\Delta_{\Shv}]
\morphism(150,615)|a|/->/<700,300>[`;\Delta_{\inf}]
\morphism(175,555)|b|/<-/<700,300>[`;\Theta_{\inf}]
\morphism(200,950)|b|/<-/<600,0>[`;\Theta_{\Shv}]
\morphism(300,50)|a|/<-/<440,0>[`;\Theta_{\Shv}]
\morphism(130,450)|b|/<-/<650,0>[`;\Theta_{\pre}]
\morphism(50,600)<0,400>[`;a\sub{\Shv}]
\morphism(1000,600)<0,400>[`;a\sub{\rel}]
\morphism(50,105)<0,400>[`;\Lambda]
\morphism(1000,100)<0,400>[`;\Phi]
\put(50,81.5){\tiny $\bot$}
\put(55,54){\tiny $\sim$}
\put(113,85){{\tiny $\vdash$}}
\put(0.5,87){{\tiny $\vdash$}}
\put(55,113.5){\tiny $\sim$}
\put(55,2){{\tiny $\sim$}}
\put(0.5,30){{\tiny $\sim$}}
\put(111.5,30){{\tiny $\sim$}}
 \efig$$

\section{ Heyting Algebras}

A complete Heyting algebra can be defined in two different ways categorically; as a partial order and as an enriched category. In this section we explore the basic concepts that underly the partial order definition and how this structure is used to define sheaves for a Hayting algebra.   The theory of presheaves and sheaves for a complete Heyting algebra is well understood (see for example \cite{mc2}\cite{Jst}), but we will go over in detail some aspects of the theory that we will borrow or mimic when we start working with relational-presheaves and relational-sheaves.  Those who are familiar with these constructions can skip this section and the next (enriched categories) and jump to section 4. The perspective  of this section is that a Heyting algebra is a partial order with additional structure.  These will be built from scratch with the focus on the structural properties that will be used later. We then move onto an exploration of the basics of sheaves for a Heyting algebra with a proof of the comparison theorem which says that for a Heyting algebra $\h$, the category of sheaves on the Heyting algebra of down-closed subsets $\D\h$, is equivalent to the category of presheaves for $\h$.

\begin{df}\em A {\it Partial order} is a category $\C$ in which the hom sets have at most one element.  Equivalently a partial order is a
pair $\langle \C, R\rangle$ where $\C$ is a set and $R$ is a reflexive and transitive relation (we don't necessarily require antisymmetry). 
\edf

Underlying the relationship between the structures of sheaves and presheaves for a Heyting algebra $\h$ is the association to a presheaf of down-closed subsets of 
$\h$ and the association to a sheaf of principal down-closed subsets. 

\begin{df} \em Let $X$ be a subset of a partial order $\ord$.

\vsc{-.5}

\begin{itemize}

\item $X$ is a {\it down-closed subset} if whenever $a\leq b\in X$, then $a\in X$.

\item The {\it down-closure} of $X$ is the set $X^{\dagger} = \{b\ |\ \exists a\in X\hbox{ and } b\leq a\ \}$

\item $X$ is a {\it principal} down-closed subset if $X = \{a\}^{\dagger}$ for some element $a\in\ord$.

\end{itemize}

\vsc{-1}

\edf

We now turn our attention to defining Heyting algebras as lattices with extra structure.

\begin{df}\em Let $\langle \L , \leq \rangle$ be a partially ordered set 

\begin{itemize}
\item $\langle \L , \leq \rangle$ is a {\it lattice} if it has all binary suprema and infima.

\item A lattice $\langle \L , \leq \rangle$ has a {\it top element} $\top$, if $\top\wedge x = x$ for all $x\in\L$.

\item A lattice $\langle \L , \leq \rangle$ has a {\it bottom element} $\bot$, if $\bot\vee x = x$ for all $x\in\L$.

\end{itemize}

\vsc{-1}

\edf

\begin{df} \em A {\it Heyting algebra}, $\h$, is a lattice with both a top and bottom element and for every $x\in\h$ the functor $()\wedge x:\h\ra\h$ has a right adjoint $x\Rightarrow ():\h\ra\h$ (called implication).  $\h$ is {\it complete} if $\h$ contains all suprema.
\edf

From this point forward when we refer to a Heyting algebra we will assume that it is complete.  There is one theorem where we explicitly use the fact that the meet operation has a right adjoint. In particular we have that $y\wedge x \leq z$ if and only if  $y\leq (x\Rightarrow z)$ and because of the symmetry of meet, $x\wedge ()$ will also have a right adjoint given by $x\wedge y\leq z$ if and only if  $y\leq (z\Leftarrow x)$.  In addition $x\!\!\iff\!\! z$ will be used to represent 
$(x\Rightarrow z)\wedge (x\Leftarrow z)$.  Since we assume a Heyting algebra has all suprema the implication operator can be determined by,

\vsc{-.5}

$$x\Rightarrow z = \bigvee\{y\ |\ y\wedge x \leq z \}.$$

The meet is symmetric and is a left adjoint so it must preserve all suprema in each variable. 

\vsc{-.75}

$$x \wedge \bigvee_{i\in I}y_i = \bigvee_{i\in I}(x\wedge y_i)\quad\quad\hbox{ and }\quad\quad \bigvee_{i\in I}x_i\wedge y = \bigvee_{i\in I}(x_i\wedge y)$$

The set of down-closed subsets of a partial order $\ord$ is a partial order.  It is easy to see that by using the intersection and union operators $\D\ord$ is a Heyting algebra.  In addition for a Heyting algebra $\h$, there is an adjunction $\bigvee\dashv()^{\dagger}:\h\ra\D\h$, where $\bigvee$ takes the supremum of a down-closed subset and $()^{\dagger}$ turns an element $a$ into its associated principal down-closed subset. This adjunction will be used explicitly when we construct the associated sheaf and  relational-sheaf functors.




\vsc{2}




\subsection{ Presheaves and Sheaves for a Heyting algebra}

For a Heyting algebra $\h$ a functor of the form $F:\h^{op}\ra\set$ is called a {\it presheaf} and if $k\leq h$ in $\h$, then there is a function $F(k\leq h)$ that sends any element in $x\in F(h)$ to an element in $F(k\leq h)(x) \in F(k)$.  $F(k\leq h)(x)$ is called the restriction of $x$ to $k$ and will be denoted $x_{|_k}$.

\begin{df}\em Let $\h$ be a Heyting algebra and  $F:\h^{op}\ra \set$ is a presheaf, then

\begin{itemize}

\item If $h\in\h$, then a {\it cover} of $h$ is a subset $A\subseteq\h$ such that $\ds\bigvee_{k\in A}\!\! k = h$

\item  A {\it matching family} for a subset $A\subseteq\h$ and a presheaf $F$ is a family of elements $\langle x_k\rangle_{k\in A}$ such $x_k\in F(k)$ for all $k\in A$
and for every $k,l \in A$ the restrictions  $x_k{_{|_{k\wedge l}}}$ and  $x_l{_{|_{k\wedge l}}}$ must be equal.

\item An {\it amalgamation} for a matching family is an element $x\in F\Big({\vee A}\Big)$ such that for every $k\in A$, the restriction of $x$ to $k$,  $x_{|_k}$, is equal to $x_k$.

\item A presheaf $F$ is a {\it sheaf} if every matching family has a unique amalgamation.

\end{itemize}

\vsc{-.75}

\edf

The category of presheaves and transformations between them will be denoted $\Preshv$ and the full subcategory of sheaves is denoted $\shv$.

Let $F:\D\h^{op}\ra\set$ be a presheaf on the Heyting algebra of down-closed subsets of $h$ and assume $x\in F(A)$, where $A$ is a down-closed subset of $\h$. Associate to $x$ the matching family  $x^{\dagger} = \{x_{|_{h^{\dagger}}}\ |\ h\in A\}$.

It is important to point out here that the  construction of $x^{\dagger}$ only uses principal down-closed subsets of $A$.   In this way we are in essence associating to the element $x$ a down-closed subset of $\h$ and a matching family.   In the comparison theorem we will utilize this and the fact that any down-closed subset of $\h$ can be represented as the union of all the principal down-closed subsets contained in it.




\begin{lm}\em Let $F:\D\h^{op}\ra\set$ be a sheaf, $x,y\in F(A)$ and $x^{\dagger} = y^{\dagger}$, then $x = y$.
\end{lm}
\pf This follows automatically since $x^{\dagger}$ and $y^{\dagger}$ are matching families for $A$ and $F$ is a sheaf.\qed

\begin{thm}(Comparison Theorem)\em\ Let $\h$ be a Heyting algebra, then the category of presheaves for $\h$ is equivalent to the category of sheaves on $\D\h$.
\end{thm}

\pf 
We begin by defining the functor $\Gamma:\Preshv\ra Shv(\D\h)$.  Let $A$ be a subset of $\h$ and $B$ any subset of $A$

\begin{itemize}
\item  {\bf On objects}: $\Gamma(F)(A) = \{X\ |\ X\hbox{ is a matching family for }A\ \}$.

\item  {\bf On arrows}: $\Gamma(F)(B\subseteq A)(X) = X_{|_B}$ where $$X_{|_B} = \{ x\in X\ |\ x\in F(k)\hbox{ and } k\in B\ \}$$
\end{itemize}

Clearly $\Gamma(F)$ is a presheaf on the down-closed subsets of $\h$.  For a transformation of presheaves $\tau:F\Rightarrow G$, then for $X$ a matching family of $A$, $\Gamma(\tau)_A$ maps $X$ to the set $\{ \tau_h(x) \ |\ x_h\in X\hbox{ and } h\in A\}$. By tracing through the appropriate diagrams element by element it is clear 
that $\Gamma(\tau)_A(X)$ is a matching family for $A$ with respect to the presheaf $G$ and that $\Gamma(\tau)$ is a transformation. For  $\twotran{F}{\tau}{G}{\sigma}{H}$ composable transformations of presheaves we have $\Gamma(\sigma\tau) = \Gamma(\sigma)\Gamma(\tau)$.

We just need to check that $\Gamma$ is a sheaf for $\D\h$.  Let $A_{i\in I}$ be a cover for $A$ ($\bigcup A_i = A$) and $X_{i\in I}$ a matching family for the cover.  Note, by construction of $\Gamma$ each $X_i$ is itself a matching family for $A_i$ and where $X_i$ and $X_j$ overlap, they must contain the same elements.  Since the supremum in $\D\h$ is simply the union of sets it easily follows that $\X = \bigcup_{\i\in I}X_i$ is a matching family for $A$ and that it is unique and that 
$\X_{|_{i}} = X_{i}$ .  Thus $\X$ is the unique amalgamation in $\Gamma(F)(A)$ for the matching family $X_i$.

Now we turn our attention to the functor $\Lambda :Shv(\D\h)\ra \Preshv\!\!.$  Let $h\in\h$ and $k\leq h$.

\begin{itemize}
\item {\bf On objects}:  $\Lambda(F)(h) = F(h^{\dagger})$

\item {\bf On arrows}: $\Lambda(F)(k\leq h) = F(k^{\dagger}\subseteq h^{\dagger})$.
\end{itemize}

Let $\tau:F\ra G$ be a transformation of sheaves then we define $\Lambda(\tau)_h$ to be the morphism $\tau_{h^{\dagger}}$.  By definition it is evident that 
$\Lambda$ is functor.

For the equivalence $\Gamma\Lambda(F) \equiv F$ we will construct for any sheaf $F\in\shvd$ 
and  each $A\in\D\h$ a bijection, $\tau_A$,  between $F(A)$ and  $\Gamma\Lambda(F)(A)$.  For each $x\in F(A)$ we set $\tau_A(x) = x^{\dagger}$.  It is straightforward to show that $\tau$ is a transformation.  Now we look at $\tau^{-1}$.
Let $X$ be a matching family for $A$ in $\Gamma\Lambda(F)(A)$.  So each $X$ has the form  $X = \langle x_h \rangle_{h\in A}$ where $x_h \in \Lambda(F)(h) = F(h^{\dagger})$.   Thus $X$ is a matching family for the set $\{ h^{\dagger}\ |\ h\in A\ \}$ which is a cover of $A$ in $\D\h$.  Since $F$ is a sheaf there is a unique amalgamation $x\in F(A)$ and $x^{\dagger} = X$.  
Using this we define $\tau_A^{-1}(X)$ to be the unique amalgamation of the matching family $X$.  $\tau_A^{-1}$ is one-to-one since if $X_1$ and $X_2$ have the same amalgamation $x$ then $X_1 =   x^{\dagger} = X_2$.  It easily follows that $\tau$ is a natural-isomorphism.

For the second equivalence $\Lambda\Gamma(F) \equiv F$, let $F$ be a presheaf for $F$ and construct for each $h\in \h$ a bijection $\sigma_h$, between $F(h)$ and $\Lambda\Gamma(F)(h)$.  Let $h\in \h$ and $x\in F(h)$ define
$\sigma_h(x) = x^{\dagger}$. 
 If $X$ is a matching family for $h^{\dagger}$, then  there must be a unique element $x\in X$ such that $x\in F(h)$ and $x^{\dagger} = X$.  We use this to define the inverse of $\sigma$;  $\sigma_h^{-1}(X) = x$ where $x$ is determined as above.  Clearly 
$\sigma$ is also a natural-isomorphism from which it follows that the category of presheaves $\Preshv$ is equivalent to the category of sheaves $\shvd$.
\qed

The comparison theorem helps to illustrate the relationship that exists between down-closed subsets and presheaves.  
In particular for a presheaf $F:\h^{op}\ra\set$ we can associate to every element $x$ in $\bigcup_h F(h)$ a down-closed subset of $\h$ and a matching family for that subset via $x^{\dagger}$ (using the obvious generalization).  When $F$ is a sheaf then for any down-closed cover $A$ of an element $h\in \h$ and an associated matching family $X_A$, the unique amalgamation is an element $x\in F(h)$ such that $x^{\dagger}$ is a matching family for $h^{\dagger}$ and $X\subseteq x^{\dagger}$.  Thus it follows that we can associate to every element in $\bigcup_h F(h)$ a principal down-closed subset of $\h$.   We thus have the property that the elements of a presheaf are down-closed subsets of $\h$ and the elements of a sheaf are principal down-closed subsets.

\section{Enriched Categories}

The definition of sheaf and presheaf for a Heyting algebra $\h$ relies on the partial order definition. It is the order-enriched interpretation of $\h$ that is utilized to define the categories of relational-presheaves and relational-sheaves.   In neither case do we utilize the fact that $\top$ is the identity element for the meet operation. So, in the enriched setting, we will only need to work with semifunctors and lax ones at that.   Here we do a quick review of the definitions and some important constructions related to order and supremum-enriched categories.  

\begin{df}\em  An {\it order-enriched} category is a category $\C$ where
\begin{itemize}

\item For every pair of objects $A,B\in\C$ the hom-set $\C(A,B)$ is a partial order.

\item For every triple of objects $A,B,C$ there is a order-preserving morphism (called composition)
$$\C_{ABC}:\C(A,B)\times\C(B,C)\ra\C(A,C)$$

\end{itemize}

\vsc{-1}

\edf

\begin{df}\em  A {\it supremum-enriched} category is a category $\C$ where
\begin{itemize}

\item For every pair of objects $A,B\in\C$ the hom-set $\C(A,B)$ is a complete lattice.

\item For every triple of objects $A,B,C$ there is a suprema-preserving morphism (called composition)
$$\C_{ABC}:\C(A,B)\otimes\C(B,C)\ra\C(A,C)$$

\end{itemize}

\vsc{-1}

\edf

Both composition morphisms are of course required to satisfy the appropriate associativity and identity conditions.
Every supremum-enriched category is obviously order-enriched.  In addition when $\V$ is a supremum-enriched category, the category $\V^{co}$ (obtained by reversing the order in each hom set) is an infimum-enriched category.  The two main enriched categories we will be working with are Heyting algebras and the category of sets and relations ($\rel$).  Both of these have an an associated involution which is the identity for a Heyting algebra and for $\rel$ it is the inverse relation.

\begin{df} \em An order-enriched category $\ord$ is {\it involutive} if there exists an order-preserving functor $()^*:\ord^{op}\ra\ord$ such that $\left(()^*\right)^* = {\mathbf  1}_{\q}$ and such a functor is called an {\it involution}.

\vsc{-.3}

\begin{itemize}
\item A morphism $f:A\ra B$ is {\it symmetric} if $f^* = f$. 

\item A morphism $f:A\ra B$ is a {\it symmetric map} if 
$$1_A \leq f^*f\quad\hbox{and}\quad ff^* \leq 1_B$$

\end{itemize}

\vsc{-.8}

\edf

The subcategory of an order-enriched category $\ord$ that has all objects and has morphisms the symmetric maps is denoted $\map^*(\ord)$


\begin{thm}\em  Let $\ord$ be an order-enriched category with involution, then $\map^*(\ord)$ is simply a category (The order on each hom-set reduces to equality).

\vsc{-.1}

\pf Assume $f\leq g:A\ra B$ and that both $f$ and $g$ are maps. Since the involution is order preserving we have that $f^*\leq g^*$.

$$g^* = 1_A\circ g^* \leq f^*\circ f\circ g^*\leq f^*\circ g\circ g^* \leq f^* \circ 1_B\leq f^*$$

Thus we also have $g^* \leq f^*$ and consequently $f = g$.
\qed

\end{thm}

\vsc{-.5}

An important construction for our work is the Karoubian envelope.  The Karou\-bian envelope is the idempotent splitting completion of a category and can be thought of as the natural way to convert a semicategory into a category.  In particular $\kar:\scat\ra\cat$ is left adjoint to the inclusion of $\cat$ into $\scat$ (Here {\bf $\scat$} represents the category semicategories).  Since the fact that $\top$ is an identity is not used for presheaves we will be using the structure of order-enriched semicategories that Heyting algebras have and define the symmetric Karoubian-envelope accordingly.

\begin{df}\em The {\it symmetric Karoubian-envelope} of an involutive order-enriched semicategory $\ord$ is the category $\kar^*(\ord)$ where 
\begin{itemize}

\item Objects: Symmetric idempotent arrows in $\C$.

\item Arrows: A morphism $\arr{f}{\phi}{g}{}$ between symmetric idempotents $\arr{A}{f}{A}{}$ and $\arr{B}{g}{B}{}$ is a morphism $\arr{A}{\phi}{B}{}$ that satisfies the two triangles 
$$\bfig
\btriangle(-250,0)|lrb|/->`->`->/<450,450>[A`A`B;f`\phi`\phi]
\put(-20,18){$=$}
\qtriangle(250,0)|alr|/->`->`->/<450,450>[A`B`B;\phi`\phi`g]
\put(60,30){$=$}
\efig$$
$$\phi\circ f = \phi = g\circ\phi.$$
\end{itemize}
\edf

We now have the tools to define the category of relational-presheaves and relational-sheaves.

\section{Relational-Presheaves and Sheaves}

Recall that to define a presheaf of a Heyting algebra $\h$ we take the point of view that $\h$ is a partial order and define a presheaf as a contravariant functor.   Our definition of relational-presheaf will take the view that a Heyting algebra is a one object supremum-enriched category where composition is given by the meet operation.  Our definition of a relational-presheaf will use $\h^{co}$ and in this way capture that a presheaf is a contravariant functor.  In addition a relational-presheaf will not care that $\h$ has an identity and so we will focus on lax-semifunctors.

With these in mind we begin our definition of  a relational-presheaf on $\h$ by looking at functions of the form $f:\h^{co}\ra\rel$.  Below we demonstrate that there is a direct relationship between the order-preserving functions $f:\h^{co}\ra\rel$ and down-closed subsets of $\h$. In addition there is a direct relationship  between those functions that preserve all infima of a subset of $\h$ and the set of principal down-closed subsets of $\h$.  This relationship is essential in how we define and how we differentiate between relational-presheaves and relational-sheaves.  

(A word on notation: A relation from set $A$ to set $B$  ($R:A\ra B$) is defined to be a subset of $B\times A$.  The composite of two relations is then defined accordingly.)

\begin{thm}\label{ordinf} \em Let $\h$ be a Heyting algebra and $A$ and $B$ sets, then a function $f:\h^{co}\ra\rel(A,B)$ is

\vsc{-.5}

\begin{itemize}

\item order-preserving, if and only if  the set $\{ h\ |\ f(h)(b,a) = 1 \}$ is a down-closed subset of $\h$ for every $a\in A$ and every $b\in B$.

\item  infima-preserving, if and only if the set $\{ h\ |\ f(h)(b,a) = 1 \}$ is a principle down-closed subset of $\h$  for every $a\in A$ and every $b\in B$.

\end{itemize}

\end{thm}

\vsc{-.2}

\pf We will prove this for infima-preserving.  

If $f$ preserves infima then we must have $f(\vee\{h\ |\ f(h)(b,a) = 1 \})(b,a) = 1$ and so the set $\{ h\ |\ f(h)(b,a) = 1\ \}$ is a principal down-closed subset of $\h$.

On the other hand if each set  $\{h\ |\ f(h)(b,a) = 1 \}$ is a principal down-closed set, then for any subset $\langle h_i\rangle$ of $\h$, we must have $f(\vee h_i) \leq \bigwedge_i f(h_i)$.  Assume that $f(h_i)(b,a) = 1$ for each $i$, then, since $\{ h\ |\ f(h)(b,a) = 1 \}$ is a principal down-closed subset of $\h$, it must be the case that  $f(\vee h_i)(b,a) = 1$\qed

There is an exercise early in Maclane\cite{mc} that asks one to  show that a natural-transformation between two functors can be defined as a family of arrows indexed by the arrows that satisfies appropriate triangles as opposed to the usual objects based definition.  It is this formulation of transformation that we take as our template to define relational-presheaves.  We start by creating a semicategory of pre-transformations which has as its objects sets and morphism families of arrows.     Recall that a function $f:\h^{co}\ra\rel$  is the same thing as a family of arrows $\langle \arr{A}{f(h)}{B}{} \rangle_{h\in\h}$ and so we can use the preceding theorem as a template for determining when a family of arrows is order/infima-preserving.

\begin{df} \em Let $\h$ be a Heyting algebra.  The category of  {\it pre-transformations} consists of

\vsc{-.5}

\begin{itemize}

\item Objects; Sets

\item Morphims; A morphism $\tau:A\ra B$ consists of an $\h$ indexed family of relations 
$$\tau:A\ra B = \langle \arr{A}{\tau_h}{B}{}\rangle_{h\in\h}$$
\begin{itemize}

\item $\tau$ is {\it order-preserving}, if the set $\{h\ |\ \tau_h(b,a) = 1\ \}$ is a down-closed set for every pair $a\in A, b\in B$.  

\item $\tau$ is {\it infima-preserving}, if the set $\{h\ |\ \tau_h(b,a) = 1\ \}$ is a principal down-closed set for every pair $a\in A, b\in B$. 

\end{itemize}

\end{itemize}

  Observe that we can associate a function $f_{\tau}:\h\ra\rel$ to a pre-transformation  $\tau:A\ra B = \langle \arr{A}{\tau_h}{B}{}\rangle_{h\in\h}$  and by Theorem \ref{ordinf}, $f_{\tau}$ is order-preserving if the family is order-preserving and $f_{\tau}$ is infima-preserving if the family is infima-preserving.  
  It must be the case that if $\tau$ is an infima-preserving pre-transformation, then for every $b\in B$ and $a\in A$, $\tau_{\bot}(b,a) = 1$ . 
 Let $\twoarr{A}{\tau}{B}{\sigma}{C}{}$ be  pre-transformations.
\begin{itemize}
\item If $\tau$ and $\sigma$ are both order-preserving, then we define the composite to be 

\vsc{-1}

\begin{eqnarray*}(\sigma\tau)_h(c,a) = 1&& \hbox{if and only if there exists } k,l\hbox{ such that } h\leq k\wedge l\hbox{ and }  \\ 
&&\sigma_k\circ\tau_l(c,a) = 1\end{eqnarray*}

\item If $\tau$ and $\sigma$ are both infima-preserving, then we define the composite to be 

\vsc{-1}

\begin{eqnarray*}(\sigma\tau)_h(c,a) = 1 & & \hbox{if and only if there exists a family } \langle k_i,l_i\rangle_{i\in I} \hbox{ such that }\\
& & h\leq \vee_i(k_i\wedge l_i)\hbox{ and }  \sigma_{k_i}\circ\tau_{l_i}(c,a) = 1\hbox{ for every }i\in I
\end{eqnarray*}

\end{itemize}

\vsc{-1}

\edf

The composition of pre-transformations is associative since the meet operation preserves suprema.  For $\h$ a Heyting algebra the category of  order-preserving pre-transformations is denoted $PT_{\Ord}(\h)$ and the category of infima-preserving pre-transformations is denoted $PT_{\inf}(\h)$.  These are categories since the pre-trans\-for\-mation $\langle \tau_h = \Delta_X\rangle_{h\in\h}$ (the diagonal relation on the set $X$), is the identity pre-transformation for the set $X$.  There is an involution on the pre-transformations defined by taking the inverse of each relation in the family  
( $\tau^{*} = \langle \arr{B}{\tau^{-1}_h}{A}{} \rangle_{h\in\h}$).  The category $PT_{\Ord}(\h)$ is an order-enriched category where $\tau\leq\sigma$, if $\tau_h\subseteq\sigma_h$ for every $h\in\h$.  In addition $PT_{\inf}(\h)$ is a supremum-enriched category where $\Big(\bigvee_{i\in I}\tau_{i}\Big)_h(b,a)  = 1$, whenever there exists $i\in I$ such that $\tau_{i_h}(b,a) = 1$.

As an aside, the category $PT_{\inf}(\h)$ is equivalent to the category of matrices on $\h$.  Recall that a matrix is a function $M:Y\times X\ra \h$.  From this we can define a pre-transformation $\tau_M:X\ra Y$ by setting for each $h\in\h$, $\tau_{M_h}(b,a) = 1$ if $h\leq M(b,a)$.  The details are easy to check (see for example \cite{wdg2}).  It follows that the construction of infima-preserving pre-transformations (which can be generalized to any quataloid\cite{pt}) is the coproduct completion of a quantaloid $\q$ as a supremum-enriched category. 



Now we start to mimic the relationships that exist between presheaves and sheaves by first proving a comparison theorem for the categories $PT_{\Ord}(\h)$ and $PT_{\inf}(\D\h)$

\begin{thm}(comparison theorem)\ \em Let $\h$ be a Heyting algebra, then $PT_{\Ord}(\h)$ is isomorphic to $PT_{\inf}(\D\h)$.
\end{thm}
\pf  First let us define $\Psi:PT_{\Ord}(\h)\lra PT_{\inf}(\D\h)$ as follows.

On any set $A$, $\Psi(A) = A$.  For a pre-transformation $\arr{A}{\tau}{B}{}$ in $PT_{\Ord}(\h)$ define for $X$ a down-closed subset of $\h$,
$\Psi(\tau)_X(b,a) = 1$ if and only if $\tau_h(b,a) = 1$ for every $h\in X$. We claim that $\Psi(\tau)$ is infima-preserving.

Let $\X = \{ X\ |\ \Psi(\tau)_X(b,a) = 1\}$ and let $Y\subseteq \bigcup\X$.  Let $h\in Y$, then there exists $X_h\in\X$ such that $h\in X_h$.  This immediately implies that $\tau_h(b,a) = 1$.  Since this is true for every such $h$, it now implies that $Y\in\X$, thus $\bigcup\X$ is a principal down-closed subset of $\h$.

Now we check that $\Psi$ is a functor.  To that end let $\sigma$ and $\tau$ be composable pre-transformations and examine $\Psi(\sigma\circ\tau)_X(c,a)$.

\begin{eqnarray*}
\Psi\big(\sigma\circ\tau\big)_X(c,a) = 1 & \hbox{ iff } & \big(\sigma\circ\tau\big)_h(c,a) =1 \hbox{ for every }h\in X\\
& \hbox{ iff } & \forall h\in X\ \exists k,l\in\h \hbox{ such that } h\leq k\wedge l\\
&&\hbox{ and } \sigma_k\tau_l(c,a) = 1\\
& \hbox{ iff } &  \forall h\in X\ \exists k,l\in\h \hbox{ such that } h\leq k\wedge l\\
&& \hbox{ and }\Psi(\sigma)_{k^{\dagger}}\Psi(\tau)_{l^{\dagger}}(c,a) = 1\\
& \hbox{ iff } & \exists\langle k_h^{\dagger},l_h^{\dagger}\rangle \hbox{ with } X \subseteq \bigcup\big(k_h^{\dagger}\wedge l_h^{\dagger}\big)\\
&& \hbox{ and } \big(\Psi(\sigma)\Psi(\tau)\big)_X(c,a) = 1
\end{eqnarray*}

Thus $\Psi\big(\sigma\circ\tau\big)$ equals $ \Psi(\sigma)\circ\Psi(\tau)$.

In reverse we have the functor $\Phi:PT_{\inf}(\D\h)\ra PT_{\Ord}(\h)$ which is defined to be the identity on objects (as $\Psi$ was).  For an infima-preserving pre-transformation $\arr{A}{\tau}{B}{}$ we have $\Phi(\tau)_h(b,a) = 1$ if and only if $\tau_{h^{\dagger}}(b,a) = 1$.  For every pair $(b,a)$ the set $\{ h \ |\ \Phi(\tau)(b,a) = 1\ \}$ is a down-closed subset of $\h$ since the set $\cup\X$ as described above is a principle down-closed subset of $\D\h$.  Thus $\Phi(\tau)$ is an order-preserving pre-transformation.  Let $\tau$ and $\sigma$ be composable pre-transformations.

\begin{eqnarray*}
\Phi\big(\sigma\circ\tau\big)_h(c,a) = 1 & \hbox{ iff } & (\sigma\circ\tau\big)_{h^{\dagger}}(c,a) = 1\\
& \hbox{ iff } &\exists \langle X_i,Y_i\rangle\in \D\h \hbox{ such that } h^{\dagger}\leq \bigcup(X_i\wedge Y_i)\\
&& \hbox{ and } \big(\sigma_{X_{i}}\circ\tau_{Y_{i}}\big)(c,a) = 1\hbox{ for every i}\\
& \hbox{ iff } &\exists i\ \exists k_i\in X_i,\exists l_i \in Y_i \hbox{ such that } h\leq k_i\wedge l_i\\
&& \hbox{ and } \Phi(\sigma)_{k_i}\Phi(\tau)_{l_i}(c,a) = 1\\
& \hbox{ iff } &  \big(\Phi(\sigma)\Phi(\tau)\big)_h(c,a) = 1
\end{eqnarray*}

Thus $\Phi$ is a functor and with simple computation we have that the composites are the appropriate identity semifunctors ($\Phi\Psi = 1$ and 
$\Psi\Phi = 1$).\qed

The ease of proving the comparison theorem here is the main reason we started with pre-transformations instead of beginning directly with relational-presheaves.  The comparison theorem for relational-presheaves will now be a simple consequence of our contructions.

\begin{df}\em Let $\h$ be a Heyting algebra, then
\begin{itemize}
\item The category of {\it relational-presheaves} for $\h$ is the category $$\Relpreshv = \map^*\Big(\skar\big(PT_{\Ord}(\h)\big)\Big)$$

\item The category of {\it relational-sheaves} for $\h$ is the category $$\Relshv = \map^*\Big(\skar\big(PT_{\inf}(\h)\big)\Big)$$
\end{itemize}

\vsc{-1}

\edf

\begin{ex}\em For the Heyting algebra $\h = (\top, \bot)$ the category $\kar^*(PT_{\inf}(\h))$ is equivalent to the category of sets and relations.  Let $F:A\ra A$ be a relational sheaf, then we can associate to it a set $X_F\subseteq A$ where $X_F$ is the set of equivalence classes determined by the partial equivalence relation $F_{\top}$ (eliminate all elements $x\in A$ where $F_{\top}(x,x) = 0$).  Given any set $X$ we can construct a relational-sheaf $G_X$ by $G_{X_{\top}}(x,y) = 1$ if and only if $x = y$.  These constructs  are  functorial (they easily extend to morphisms) and it forms an equivalence via the familiy of  relations determined by $R_F:F\ra G_{X_F}$ where $R_{\top} = F_{\top}$.   It follows that the category of relational-sheaves $\rel^{(\top,\bot)}_{\inf}$ is equivalent to the category of sets and functions.
\eex

\begin{cor}(comparison theorem) \em Let $\h$ be a Heyting algebra, then the category of relational-presheaves for $\h$ is equivalent to the category of relational-sheaves for $\D\h$.
\end{cor}

A relational-presheaf $\tau$ is thus a family of morphisms $\langle \tau_h\in\rel(A,A)\rangle$, for some set $A$, that satisfies the following conditions

\vsc{-.5}

\begin{enumerate}

 \item For every pair $k,l\in\h$\quad $\tau_l\circ\tau_k \leq \tau_{l\meet k}$.
 
 \item For every $h\in\h$, $\tau = \tau^*$.
 
 \item $\tau\circ\tau = \tau$

\end{enumerate}

Property 1 says that $\tau$ is a lax-semifunctor. The second condition says it is symmetric and the third says it is an idempotent pre-transformation.   A morphism of relational-presheaves $\theta:\tau\Rightarrow \sigma$ is a family of arrows $\langle \theta_h\in \rel(B,A)\rangle$ where $A$ is the set associated to $\tau$ and $B$ is the set associated to $\sigma$ and 

\begin{enumerate}

\item For every pair $k,l\in\h$\quad $\begin{array}{rcl}\theta_l\circ\tau_k&\leq& \theta_{l\meet k}\\
											\sigma_l\circ\theta_k&\leq&\theta_{l\meet k}\end{array}$ 

\item $\theta\circ \tau = \theta = \sigma\circ \theta$

\item $\tau \leq \theta^*\circ\theta$ and $\theta\circ\theta^*\leq \sigma$.

\end{enumerate}

Using the family of arrows version of a transformation these say that $\theta$ is a lax-transformation for which $\tau$ and $\sigma$ are identities and it is also a symmetric map.  Before we explicitly define these the main idea is the following.

\begin{thm}\em Let $\h$ be a Heyting algebra, then

\begin{itemize}

\item The category of relational-presheaves is equivalent to the category of symmetric idempotent order-preserving lax-semifunctors of the form $F:\h^{co}\ra\rel$ and suitable lax-transformations $\tau:F\Rightarrow G$.

\item The category of relational-sheaves is equivalent to the category of symmetric idempotent infima-preserving lax-semifunctors of the form $F:\h^{co}\ra\rel$ and suitable lax-transformations $\tau:F\Rightarrow G$.

\end{itemize}

\end{thm}

\begin{df}\em  Let $\h$ be a Heyting algebra

\begin{itemize}

\item A {\it lax-semifunctor} $F:\h^{co}\lra\rel$ consists of

\begin{itemize}

\item A set which we denote $F(*)$.

\item A function $F:\h^{op}(*,*)\lra\rel$, such that for any two elements\\ $h,k\in\h$, $$F(h)\circ F(k)\leq F(h\wedge k).$$

\end{itemize}

\item If $F$ is a lax-semifunctor, then

\begin{itemize}

\item $F$ is order-preserving, if whenever $h\leq k$, then $F(k)\leq F(h)$.

\item $F$ is infima-preserving, if $F(\vee h_i) = \wedge F(h_i)$

\item $F$ is {\it symmetric}, if for every $h\in\h$, $F(h) = F(h)^{-1}$.

\end{itemize}

\item A {\it lax-transformation}, $\tran{F}{\tau}{G}$,  consists of a family of relations, $$\langle \tau_h : K(*)\ra L(*) \rangle_{_{h\in\h}}$$ indexed by the elements of $\h$ such that for any element $h, k$ in $\h$  the following triangles hold.

\vsc{-.75}

$$\bfig
\btriangle(-350,0)|lrb|/->`->`->/<450,450>[F(*)`F(*)`G(*);F(h)`\tau_{_{h\wedge k}}`\tau_{_k}]
\put(-30,18){$\subseteq$}
\qtriangle(350,0)|alr|/->`->`->/<450,450>[F(*)`G(*)`G(*);\tau_{_h}`\tau_{_{h\wedge k}}`G(k)]
\put(70,30){$\supseteq$}
\efig$$

\begin{itemize}

\item We say that $\tau$ is order-preserving if for every $a\in F(*)$ and $b\in F(*)$ the set $\{ h\ |\ \tau_h(b,a) = 1\ \}$ is a down-closed set.

\item We say that $\tau$ is infima-preserving if for every $a\in F(*)$ and $b\in F(*)$ the set $\{ h\ |\ \tau_h(b,a) = 1\ \}$ is a principal down-closed set. 
\end{itemize}
\end{itemize}

\vsc{-1}

\edf

Let $\twotran{F}{\tau}{G}{\sigma}{H}$ be composable lax-transformations, then the composite $\sigma\circ\tau$ is defined depending on whether we wish to preserve the order or the infima as follows

\begin{itemize}

\item Order Preserving: If $\tau$ and $\sigma$ are order preserving and if we wish the composite to also be order preserving, then we define $\sigma\tau$ by
$(\sigma\tau)_h(c,a) = 1$ if and only if there exists $k,l\in\h$ such that $l\leq k\wedge l$ and $\sigma_k\tau_l(c,a) = 1$.

\item Infima Preserving:  If $\tau$ and $\sigma$ are infima-preserving and if we wish the composite to also be infima-preserving, then we define $\sigma\tau$ by
$(\sigma\tau)_h(c,a) = 1$ if and only if there exists a family of morphisms $\langle k_i,l_i\rangle$ from $\h$ such that $h\leq \vee(k_i\wedge l_i)$ and $\sigma_{k_i}\tau_{l_i}(c,a) = 1$ for every $i$.

\end{itemize}
Any lax-semifunctor $F$ is a lax-transformation $\tau_F:F\ra F$, where $\tau_F$ is the family of morphisms determined by $\langle F(h)=\tau_{F_h}\rangle_{h\in\h}$.  Since $F$ is a lax-semifunctor it is automatic that the two defining triangles for a lax-transformation hold and they both represent  $\tau_F\circ\tau_F\leq\tau_F$.  But notice that $\tau_F$ need not be an identity for $F$.  We denote the related semicategory by $\lax\rel_{\cal V}{(\h)^{co}}$.  Let $F$ and $G$ be lax-semifunctors, then 
the hom set $\lax\rel_{\cal V}^{(\h)^{co}}\!\!(F,G)$ is a partial order where $\tau\leq\sigma$ if $\tau_h\leq\sigma_h$ for every $h\in\h$.  There is also an involution on the transformations defined by $\tau^{\circ} = \langle\tau^{-1}\rangle_{h\in\h}$.  It is straight forward to show that

\vsc{-.25}

\begin{itemize}
\item[] \hspace*{-.95cm}{$\bullet$}\ The category of relational-presheaves, $\Relpreshv$, is $\map^*(\kar^*(\lax\rel_{\Ord}(\h)^{co}))$.

\item[] \hspace*{-.95cm}{$\bullet$}\  The category of relational-sheaves, $\Relshv$, is $\map^*(\kar^*(\lax\rel_{\inf}(\h)^{co})).$
\end{itemize}

Let $F$ be a relational-sheaf, then since $F$ is infima-preserving it is automatically order-preserving and thus every relational-sheaf is an order preserving lax-semifunctor.  We have a problem in that $\iota(F)$ ($F$ interpreted as order-preserving) need not be a relatonal-presheaf because the definition of the composition changes sufficiently so that the $\iota$ need not preserve composition. The only thing that is guaranteed is that $\iota(F)\iota(F)\leq\iota(F)$. So $\iota(F)$ need not be an idempotent.

Now if $F$ is a relational-presheaf then it is possible to use $\bigvee$ to turn it into an infima-preserving lax-semifunctor by setting $\bigvee(F)(h)(b,a) = 1$ if and only if 
$h\leq\vee\{ k\ |\ F(k)(b,a) = 1\ \}$. It is not difficult to show that $\bigvee(F)\bigvee(F)\leq\bigvee(F)$, but it need not be the case that $\bigvee(F)$ is an idempotent.
So this simple construction can not be interpreted as an associated relational-sheaf functor.  But fortunately we can use it in conjunction with other constructs to create an associated relational-sheaf functor.








\subsection{Singletons}

The main tool we will be using when we come to constructing the associated relational-sheaf functor and the equivalences with the sheaf constructions is singletons and singleton morphisms.  Very basically a singleton will be a relational-sheaf that represents a single element $h\in\h$.  A singleton morphism 
is a lax-transformation that simultaneously gives us an amalgamation of a matching family for a cover of $h$ determined by a given relational-sheaf $F$.  

Let $h\in\h$, then there is a relational-sheaf $F_h$, where 
\begin{itemize}
\item $F_h(*) = \{*\}$
\item $F_h(k)(*,*) = 1$ if $k\leq h$.
\end{itemize}

\begin{df}\em Let $\h$ be a Heyting algebra, then

\begin{itemize}

\item The relational-sheaves of the form $F_h$ will be called {\it singletons}

\item Let $F$ be a relational-sheaf, then any morphism of the form $\tran{F_h}{\alpha}{F}$ is called a {\it singleton morphism}.

\end{itemize}

\vsc{-1}

\edf

The following examples and constructs using singletons and singleton morphisms capture important relationships that we will constantly utilize when constructing our series of functors comparing presheaves, relational-presheaves, sheaves and relational-sheaves.

\begin{ex}\em Let $\h$ be a Heyting algebra then

\begin{itemize}

\item Let $h\leq k$, then there is a singleton morphism $\alphasingle{h,k}{}:F_h\Ra F_k$ defined by
$\alphasingle{h,k}{l}(*,*) = 1 $ if $l\leq h$.

\item  Let $\tran{F_k}{\alpha}{F}{}$ be a singleton, then $\alpha$ restricted to $h$ is the singleton morphism 
$$\alpha_{|_h} = \twotran{F_h}{\alphasingle{h,k}{}}{F_k}{\alpha}{F}.$$

\item Let $F$ be a relational-sheaf  then for each $x\in F(*)$, there is an associated representable singleton morphism $\tran{F\sub{x}}{\alphasingle{x}{}}{F}$ where 

\begin{itemize}
\item $F\sub{x}(*) = \{ *\}$ and $F\sub{x}(l)(*,*) = 1$ if $F(l)(x,x) = 1$

\item $\alphasingle{x}{l}(y,*) = 1$ if $F(l)(y,x) = 1 $.

\end{itemize}

\item If $\alpha:F_h\ra F$ is a singleton morphism, then there is a morphism of the lax-semifunctors $\alpha^{\circ}:F\ra F_h$ determined by the involution \Big(note: it need not be a morphism in $\Relshv$\Big).  

\item If $\Tran{\alpha}{F_h}{F}$ is a singleton, then we can associate to each $x\in F(*)$ an element $\alpha(x)\in\h$ by setting $\alpha(x) = \vee\{h\ |\ \alpha_h(x,*) = 1\ \}$.  Note that $(\alpha^{\circ}\alpha)_{l}(*,*) = 1$ if and only if $l\leq\vee_y\alpha(y)$.

\end{itemize}

\vsc{-1.3}

\eex

\begin{lm}\em If $\tran{F_h}{\alpha}{F}$ is a singleton in $\Relshv$, then $\alpha$ is a monomorphism.
\end{lm}
\pf By assumption $\alpha F_h = \alpha$ and $F_h \leq \alpha^{\circ}\alpha$.  The first assumption tells us that $\alpha_l(y,*) = 1$, if $l\leq \alpha(y) = \bigvee\{k\wedge h\  |\ \alpha_k(y, *) = 1\ \}\leq h$.  We know that $F_h(h)(*,*) = 1$, which tells us that $(\alpha^{\circ}\alpha)_h(*,*) = 1$.  This implies that $h \leq \vee_y \alpha(y)$ and thus $h = \vee_y \alpha(y)$.  It follows that $(\alpha^{\circ}\alpha)_k(*,*) = 1$ if and only if $k\leq h$ and thus $F_h = \alpha^{\circ}\alpha$.\qed

Now we prove a result that is central to our construction of the associated sheaf and relational-sheaf functors.  In particular it allows us to transform back and forth between a sheaf construct (the restrictions) and the relational-sheaf setting (the singleton morphisms $\alpha^{\circ}\beta$).
\begin{lm}\label{lm4}\em Let $\tran{F_h}{\alpha}{F}{}$ and $\tran{F_k}{\beta}{F}{}$ be singletons in $\Relshv$, then $\alpha_{|_l} = \beta_{|_l}$ if and only if 
$(\alpha^{\circ}\beta)\sub{l}(*,*) = 1$.
\end{lm}

\pf  We first claim that for each $x\in F(*)$ that $\alpha_{|_l} = \beta_{|_l}$, if and only if $$l\leq h\wedge k \wedge \bigwedge_x(\alpha(x) \iff \beta(x))$$

\begin{eqnarray*}
\alpha_l  =  \beta_l & \hbox{ iff } &  \forall x,\quad l \leq h \wedge k \hbox{ and } \\
& &  \quad\quad\quad\quad(\alpha\circ\alphasingle{l,h}{})_l(x,*) = (\beta\circ\alphasingle{l,k}{})_l(x,*),\\
& \hbox{ iff } & \forall x,\quad l\leq h\wedge k \hbox{ and } \alpha(x)\wedge l = \beta(x) \wedge l\\
& \hbox{ iff } & \forall x,\quad l\leq h\wedge k \hbox{ and } l\wedge \alpha(x) \leq \beta(x)\hbox{ and }l\wedge\beta(x) \leq \alpha(x)\\
&\hbox{ iff }& \forall x,\quad l\leq h\wedge k \hbox{ and } l\leq (\alpha(x)\iff\beta(x))\\
& \hbox{iff} & \ \ \quad\quad l\leq h\wedge k\wedge \bigwedge_x(\alpha(x)\iff\beta(x)) 
\end{eqnarray*}

Let us recall that $(\alpha^{\circ}\beta)_l(*,*) = 1$ if and only if $l\leq \vee\{r\wedge s\ |\ \alpha^{\circ}\sub{r}\beta\sub{s}(*,*) = 1\}$.  But this is equivalent to asking that $l\leq \vee_x\{\alpha(x)\wedge\beta(x)\}$.  We now claim that $ l\leq h\wedge k\wedge \bigwedge_x(\alpha(x)\iff\beta(x)) $ if and only if $l\leq \vee\{\alpha(x)\wedge\beta(x)\}$ from which the lemma follows.

Clearly $\Big(h\wedge k \wedge \bigwedge_x\alpha(x)\iff\beta(x)\Big) \leq \bigvee_x\{\alpha(x)\wedge\beta(x)\}$.  For the inequality in the other direction observe that since $\alpha$ and $\beta$ are maps we know that $\alpha(x)\leq h$ and $\beta(x)\leq k$.  For any $x\in F(*)$ we have the following.

\begin{eqnarray*}
(\alpha^{\circ}\beta\beta^{\circ})\sub{l}(*,x) = 1 & \hbox{iff} & l\leq\bigvee_y \{r\wedge s\wedge t\ |\ \alpha_r^{\circ}(*, y)\wedge\beta_s(y,*)\wedge\beta_t^{\circ}(*,x) = 1\ \}\\
&\hbox{ iff } & l\leq \bigvee_y \alpha(y)\wedge\beta(y)\wedge\beta(x)\\
\end{eqnarray*}

Since $\beta$ is a map we have $\alpha^{\circ}\beta\beta^{\circ}\leq\alpha^\circ F\leq \alpha^{\circ}$.  Thus
$$ l\leq \bigvee_y \alpha(y)\wedge\beta(y)\wedge\beta(x)\leq\alpha(x)$$
which is equivalent to
$$ l\leq \bigvee_y \alpha(y)\wedge\beta(y)\leq\beta(x)\Rightarrow \alpha(x)$$

Similarly $l\leq  \bigvee_y\alpha(y)\wedge\beta(y) \leq \beta(x)\Ra \alpha(x)$.  These are true for any $x$ so we must have that 
$$l\leq \bigvee_x\{\alpha(x)\wedge\beta(x)\} \leq h\wedge k\wedge \bigwedge_x\Big(\alpha(x)\iff\beta(x)\Big)$$
 From which it now follows that $\alpha_{|_l} = \beta_{_l}$ if and only if $(\alpha^\circ\beta)\sub{l}(*,*) = 1$\qed


\section{The Adjunction $\Delta_{\inf}\dashv\Theta_{\inf}:\Relshv\ra\Preshv$.}

















Our goal is to create the equivalences from which it follows that presheaves and sheaves can be interpreted as being lax-semifunctors with added structure.  We begin by focusing on constructing an adjunction between the category of presheaves and the category of relational-sheaves $(\Delta_{\inf}\dashv\Theta_{\inf})$.


\begin{df}\em  Let $F:\h^{co}\ra\rel$ be a relational-sheaf and $\tran{F}{\tau}{G}$ a morphism of relational-sheaves, then the functor 
$\Theta_{\inf}:\Relshv\ra \Preshv$ is defined as follows

\begin{itemize}

\item {\bf On Objects:}

\begin{itemize}

\item $\Theta_{\inf}(F)(h) = \{ \tran{F_{h}}{\alpha}{F}\ |\ \alpha \hbox{ is a singleton }\ \}$.
\item If $k\leq h$, then $\Theta_{\inf}(F)(k\leq h)(\alpha) = \Big(\twotran{F_{k}}{\alphasingle{k,h}{}}{F_{h}}{\alpha}{F}\Big) = \alpha_{|_k}$.
\end{itemize}

\item {\bf On Arrows:}  $$\Theta_{\inf}(\tau)_h(\alpha) = \Big(\twotran{F_{h}}{\alpha}{F}{\tau}{G}\Big) $$

\end{itemize}

\vsc{-1}

\edf


 Clearly by definition, $\Theta_{\inf}(F)$ is a presheaf and $\Theta_{\inf}(\tau)_h$ is a function for each $h\in\h$.  In addition since both $\Theta_{\inf}(F)(h\leq k)$ and $\Theta_{\inf}(\tau)_h$ are defined as the composition of morphisms it is straight forward to show that $\Theta_{\inf}(\tau)$ is a transformation and for composable morphisms that $\Theta_{\inf}(\sigma\tau) = \Theta_{\inf}(\sigma)\circ\Theta_{\inf}(\tau)$.  

Now we turn our attention to $\Delta_{\inf}$   

















\begin{df}\em Let $F:\h^{op}\ra\set$ be a presheaf and $\tran{F}{\tau}{G}$ be a transformation of presheaves, then the functor $\Delta_{\inf}:\Preshv\!\!\ra\Relshv$ is defined as follows

\begin{itemize}

\item {\bf On Objects:}

\begin{itemize}

\item $\ds\Delta_{\inf}\Big(F\Big)(*) = \coprod_{h\in\h}F(h)$

\item $\ds\Delta_{\inf}\Big(F\Big)(h)(b,a) = 1$ if and only if $$h\leq\bigvee\{k\wedge l\ |\ b\in F(l), a\in F(k) \hbox{ and } b_{|_{k\wedge l}} = a_{|_{k\wedge l}}\Big\}$$

\end{itemize}

\item {\bf On Morphisms:}

\begin{itemize}

\item $\Delta_{\inf}\big(\tau\big)_h(b,a) = 1$ if and only if 
$$\ds h\leq \bigvee\Big\{k\wedge l\ |\ b\in G(l), a\in F(k) \hbox{ and } b_{|_{k\wedge l}} = \tau_k(a)_{|_{k\wedge l}}\Big\}$$

\end{itemize}

\end{itemize}

\vsc{-1}

\edf

A particularly useful relationship is that  $\Delta_{\inf}(G)(h)(b,\tau_k(a)) = 1$ if and only if $\Delta_{\inf}(\tau)_h(b,a) = 1$.  This is automatic since both rely on the fact that $b_{|_{k\wedge l}} = \tau_k(a)_{|_{k\wedge l}}$.

To check the details that $\Delta_{\inf}$ as defined is a functor we first observe that by construction $\Delta_{\inf}(F)$ is clearly symmetric and infima-preserving and $\Delta_{\inf}(\tau)$ is infima-preserving as well.   First we show that $\Delta_{\inf}(\sigma\tau) = \Delta_{\inf}(\sigma)\Delta_{\inf}(\tau)$ for composable transformations
 $\twotran{F}{\tau}{G}{\sigma}{H}$.

$\Delta_{\inf}(\sigma\tau)_h(c,a) = 1$ if and only if 
$$\ds h\leq \bigvee\Big\{k\wedge l\ |\ c_{|_k}\in H(k), a_{|_l}\in F(l) \hbox{ and } c_{|_{k\wedge l}} = (\sigma\circ\tau)_{k\wedge l}(a)_{|_{k\wedge l}}\Big\}$$

Now $\Big(\Delta_{\inf}(\sigma)\Delta_{\inf}(\tau)\Big)_h(c,a) = 1$ if and only if there exists a family of morphisms $\langle k_i,l_i\rangle_{i\in I} $ such that $h\leq  \vee_i(k_i\wedge l_i)$
and
$\Delta_{\inf}(\sigma)_{k_i}\Delta_{\inf}(\tau)_{l_i}(c,a) = 1$.  


This is if and only if  there exists a family of elements $b_i\in G(r_i)$  such that  $\Delta_{\inf}(\sigma)_{k_i}(c,b_i) = 1$ and $\Delta_{\inf}(\tau)_{l_i}(b_i,a) = 1$.

{\tiny
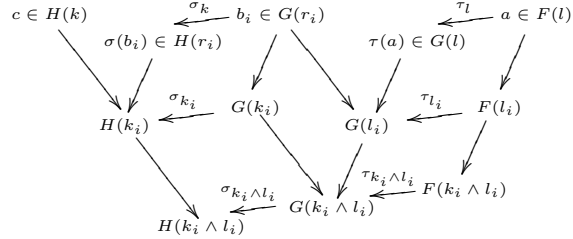
$$\bfig
\put(-285,0){$c\in H(k)$}
\put(-250,-10){$\sigma(b_i)\in H(r_i)$}
\put(-200,0){$b_i\in G(r_i)$}
\put(-150,-10){$\tau(a) \in G(l)$}
\put(-100,0){$a\in F(l)$}
\put(-252,-42){$H(k_i)$}
\put(-202,-35){$G(k_i)$}
\put(-159,-42){$G(l_i)$}
\put(-109,-35){$F(l_i)$}
\put(-230,-80){$H(k_i\wedge l_i)$}
\put(-180,-73){$G(k_i\wedge l_i)$}
\put(-130,-66){$F(k_i\wedge l_i)$}

\morphism(-2300,0)<225,-290>[`;]
\morphism(-1950,-75)<-100,-225>[`;]
\morphism(-1700,10)<-200,-25>[`;\sigma_k]
\morphism(-1550,-0)<-100,-225>[`;]
\morphism(-1550,-0)<225,-290>[`;]
\morphism(-1150,-75)<-100,-225>[`;]
\morphism(-850,10)<-200,-25>[`;\tau_l]
\morphism(-750,-0)<-100,-225>[`;]
\morphism(-1750,-300)<-200,-25>[`;\sigma_{k_i}]
\morphism(-950,-300)<-200,-25>[`;\tau_{l_i}]
\morphism(-2050,-350)<225,-290>[`;]
\morphism(-1650,-275)<225,-290>[`;]
\morphism(-1550,-600)<-170,-20>[`;\sigma_{k_i\wedge l_i}]
\morphism(-1275,-350)<-100,-225>[`;]
\morphism(-875,-275)<-100,-225>[`;]
\morphism(-1100,-550)<-170,-20>[`;\tau_{k_i\wedge l_i}]
\efig$$

}

 All we need to do is to let each $b_i$ be $\tau_{l}(a)$ and then, by the diagram above, this reduces to finding $\langle k_i,l_i\rangle\in\h$ with $h\leq \vee_i (k_i\wedge l_i)$ and $c_{|_{{k_i}\wedge {l_i}}} = \sigma_k(\tau_k(a))_{|_{{k_i}\wedge {l_i}}}$ for all $i\in I$. But this is equivalent to our requirement for $\Delta_{\inf}(\sigma\tau)$.

It now follows that $\Delta_{\inf}(F)$ is an idempotent and that $\Delta_{\inf}(\tau)$ is a morphism.  To see this we simply observe that a presheaf $F$ is simultaneously the identity transformation for $F$ which implies that $\Delta_{\inf}(\tau)\Delta_{\inf}(F) = \Delta_{\inf}(\tau) = \Delta_{\inf}(G)\Delta_{\inf}(\tau)$ and that  
$\Delta_{\inf}(F)$ is an idempotent.

Finally to see that $\Delta_{\inf}(\tau)$ is a map  we observe that if $\Delta_{\inf}(F)(h)(c,a) = 1$, then there exists a family $\langle k_i,l_i\rangle\in\h$ where $h\leq \vee_i (k_i\wedge l_i)$ and $c_{|_{k_i\wedge l_i}} = a_{|_{k_i\wedge l_i}}$.  Since $\tau$ is a transformation we know that $\tau_{k_i}(c)_{|_{k_i\wedge l_i}} = \tau_{l_i}(a)_{|_{k_i\wedge l_i}}$.  Thus both $\Delta_{\inf}(\tau)^{\circ}_{k_i}(c, \tau_{k_i(}a))$ and $\Delta_{\inf}(\tau)_{l_i}(\tau_{k_i(}a), a)$ are equal 1 and hence $(\Delta_{\inf}(\tau)^{\circ}\Delta_{\inf}(\tau))_h(c,a) = 1$.  Thus $\Delta_{\inf}(F) \leq \Delta_{\inf}^{\circ}(\tau)\Delta_{\inf}(\tau)$. Now the equality $\Delta_{\inf}(\tau)_h(b,a) = \Delta_{\inf}(G)(h)(b,\tau_h(a))$ can be used to show, in a similar vein, that $\Delta_{\inf}(\tau)\Delta_{\inf}(\tau)^{\circ} \leq \Delta_{\inf}(G)$.

We have shown that $\Delta_{\inf}$ and $\Theta_{\inf}$ are functors and now to show that $\Delta_{\inf}\dashv\Theta_{\inf}$.


We start with the unit of the adjunction, $\eta:{\mathbf 1}\Ra\Theta_{\inf}\Delta_{\inf}$. 
Let $F$ be a presheaf and define $\eta_F:F\Ra\Theta_{\inf}\Delta_{\inf}(F)$ by setting $\eta_{F,h}(x) = \alphasingle{x}{}$.  We need to show that both $\eta$ and $\eta_F$ are transformations as exhibited by the following squares.

$$\bfig
\square(-600,500)|ammb|/>`>`>`>/<650,500>[F`\Theta_{\inf}\Delta_{\inf}(F)`G`\Theta_{\inf}\Delta_{\inf}(G);\eta_F`\tau`\Theta_{\inf}\Delta_{\inf}(\tau)`\eta_G]
\square(600,500)|ammb|/>`>`>`>/<700,500>[F(h)`\Theta_{\inf}\Delta_{\inf}(F)(h)`F(k)`\Theta_{\inf}\Delta_{\inf}(F)(k);\eta_{F,h}`k\leq h`\Theta_{\inf}\Delta_{\inf}(F)(k\leq h)`\eta_{F,k}]
 \efig$$

Starting with $\eta_F$ we need to know that if $x\in F(h)$, then the singletons, $\Delta_{\inf}(F)_x$ and $\Delta_{\inf}(F)_h$ are equal. But  this is automatic since $\Delta_{\inf}(F)(l)(x,x) = 1$ if and only if $x_{|_{l}} = x_{|_{l}}$,  which will only occur whenever $l\leq h$.  And so the square on the right requires for every $x\in F(h)$ that the singletons,  $\twotran{\Delta_{\inf}(F)\sub{k}}{\alphasingle{k,h}{}}{\Delta_{\inf}(F)\sub{h}}{\alphasingle{x}{}}{\Delta(F)}$ and $\tran{\Delta_{\inf}(F)\sub{k}}{\alphasingle{x}{}{{|_{k}}}\ \ }{\Delta_{\inf}(F)}$ be equal.  Tracing through the details we find that both $(\alpha\sub{x}\alphasingle{k,h}{})_l(y,*) =1$ and $(\alpha\sub{x_{|_k}}){_l}(y,*) = 1$ simply require that $y_{|_l} = x_{|_l}$.  Similarly $\eta$ is a transformation because  $(\Theta_{\inf}\Delta_{\inf}(\tau)\circ\eta_F)(x) = \twotran{F_h}{\alpha_x}{\Delta_{\inf}(F)}{\Delta_{\inf}{\tau}}{\Delta_{\inf}(G)}$ and $\eta_G\circ\tau = \tran{F_h}{\alpha\sub{\tau\sub{h}\!\!(x)}}{\Delta_{\inf}(G)}$ are equal since both cases  reduce down to the requirement that  $y_{|_l} = \tau\sub{h}(x){{_{|_l}}}$.

The counit of our adjunction $\varepsilon:\Delta_{\inf}\Theta_{\inf}\ra {\mathbf 1}:\Relshv\ra\Relshv$ is defined for a relational-sheaf $F$  by letting the morphism
$\varepsilon_F:\Delta_{\inf}\Theta_{\inf}(F)\Ra F$ be given by

$$\varepsilon_{F,h}(x,\alpha) = 1\hbox{ if and only if }\alpha_h(x, *) = 1.$$

First we will show that for each $F$, $\varepsilon_{F,h}$ is not just a morphism, but an isomorphism.

$(\varepsilon_{ F}\tau\sub{\Delta\Theta F})_h(x,\alpha) = 1$ if and only if there exists a family of pairs of elements in $\h$ $\langle k_i,l_i\rangle$, such that $h\leq \vee(k_i\wedge l_i)$ and 
$\varepsilon_{k_i}\tau\sub{\Delta\Theta F_{l_i}}(x,\alpha) = 1$.  

This occurs if and only if there exists a family $\langle k_i,l_i\rangle$, a singleton morphism $\beta$ such that $h\leq \vee(k_i\wedge l_i)$ and $\varepsilon_{k_i}(x,\beta)= 1$ and $\tau\sub{\Delta\Theta F_{l_i}}(\beta,\alpha) = 1$.  But these just say that 
$\alpha_{k_i}(x,*) =1$ and $\beta_{l_i} = \alpha_{l_i}$ so picking $h = k_i = l_i$ and $\beta = \alpha$ gives the desired result that $$\varepsilon_{ F}\tau\sub{\Delta\Theta F} = \tau\sub{\Delta\Theta F}.$$ Since each $\alpha$ is a singleton morphism it follows that  $$\tau_F\varepsilon_F = \varepsilon_F.$$

$\varepsilon_F$ is a symmetric map (isomorphism) since $(\varepsilon_F^{\circ}\varepsilon_F)_h(\beta,\alpha) = 1$, if and only if there exists  a family $\langle k_i,l_i\rangle$ such that $h\leq \vee(k_i\wedge l_i)$ and 
$\varepsilon_{F_k}^{\circ}\varepsilon_{F_l}(\beta,\alpha) = 1$.  

But this happens if and only if there exist a family $\langle k_i,l_i\rangle$  and associated to each pair $(k_i, l_i)$ there is an $x_i$ such that $h\leq \vee(k_i\wedge l_i)$ and $\beta_{k_i}(x_i,*) = 1$ and $\alpha_{l_i}(x_i.*) = 1$.  

This is equivalent to saying that there exist a family $\langle k_i,l_i\rangle$  such that $h\leq \vee(k_i\wedge l_i)$ and
$\beta^{\circ}_{k_i}\alpha_{l_i}(*,*) = 1$.  By lemma \ref{lm4} this is if and only if there exist a family $\langle k_i,l_i\rangle$  such that $h\leq \vee(k_i\wedge l_i)$ and 
 $\beta_{|_{h}} = \alpha_{|_{h}}$.  Thus we have $(\varepsilon_F^{\circ}\varepsilon_F)_h(\beta,\alpha) = 1$ 
 if and only if $\Delta_{\inf}\Theta_{\inf}(F)(h)(\beta,\alpha) = 1$ and so
 $\varepsilon_F^{\circ}\varepsilon_FÊ= \Delta_{\inf}\Theta_{\inf}(F)$

 For the other inequality we observe that $(\varepsilon_F\varepsilon_F^{\circ})(h)(c,a) = 1$ if and only if there exists a family $\langle k_i,l_i\rangle$ and associated to each pair $(k_i, l_i)$ there is a singleton $\alpha^i$ such that 
$h\leq \vee(k_i\wedge l_i)$ and $\alpha^i_{k_i}\alpha^{i^{\circ}}_{l_i}(c,a) = 1$.  Since $\alpha^i$ is a singleton $\alpha^i\alpha^{i^{\circ}} \leq F$ from which we conclude that $F(h)(c,a) = 1$.  Thus $\varepsilon_F\varepsilon_F^{\circ}\leq F$

Now we observe that if $F(h)(c,a) = 1$ then $F(h)(c,c) = 1$ since $F$ is an idempotent.
 For each $\alpha^i$ pick the representable singleton $\alphasingle{c}{h}$ to obtain $\alphasingle{c}{h}\alphasingle{c^{\circ}}{h}(c,a) = 1$.   Now we have the inequality in the other direction  $F\leq \varepsilon_F\varepsilon_F^{\circ}$ thus
 
 $$F = \varepsilon_F\varepsilon_F^{\circ}$$


Finally to show that  $\varepsilon$ is a transformation, and thus a natural isomorphism, the following square must commute for each presheaf $F$.

$$\bfig
\square(-500,500)|ammb|/>`>`>`>/<700,500>[\Delta_{\inf}\Theta_{\inf}(F)`F`\Delta_{\inf}\Theta_{\inf}(G)`G;\varepsilon_F`\Theta_{\inf}\Delta_{\inf}(\tau)`\tau`\varepsilon_G]
 \efig$$

The right side of the square gives us that \noindent $(\tau\circ\varepsilon_F)_h(b,\alpha) = 1$ if and only if there exists a family $\langle k_i,l_i\rangle_{i\in I}$ such that $h\leq \vee (k_i\wedge l_i)$ and $\tau_{k_i}\circ\varepsilon_{F_{l_i}}(b,\alpha) = 1$ for each $i\in I$. Which happens if and only if there exists a family  $\langle k_i,l_i\rangle_{i\in I}$ such that $h\leq  \vee (k_i\wedge l_i)$ and $\tau_{k_i}\circ\alpha_{l_i}(b,*)= 1$,  and thus $(\tau\circ\varepsilon_F)_h = (\tau\circ\alpha)_h$.

On the left side of the square; \noindent $(\varepsilon_G\circ\Delta_{\inf}\Theta_{\inf}(\tau))_h(b,\alpha) = 1$, if and only if there exists a family $\langle k_i,l_i\rangle_{i\in I}$ such that $h\leq \vee (k_i\wedge l_i)$  and $(\varepsilon_{G_{k_i}}\circ\Delta_{\inf}\Theta_{\inf}(\tau)_{l_i})(b,\alpha) = 1$ for each $i\in I$. 

This is equivalent to the existence of a family $\langle k_i,l_i\rangle_{i\in I}$ such that $h\leq \vee (k_i\wedge l_i)$ and for each $i\in I$ a $\beta_i$ such that 
$\varepsilon_{G_{k_i}}(b,\beta_i) = 1$ and $\Delta_{\inf}\Theta_{\inf}(\tau)_{l_i}(\beta,\alpha) = 1$


Which is in turn equivalent to saying that there exists a family $\langle k_i,l_i\rangle_{i\in I}$ such that $h\leq \vee (k_i\wedge l_i)$ and associated to each pair 
$(k_i,l_i)$, there is a singleton morphism $\beta_i$ where $\beta_{i_k}(b,*) = 1$ and $\tau\alpha_{|_l} = \beta_{i_{|_l}}$.

 But of course $\tau\alpha$ is such a 
 $\beta_i$ and $h\leq \vee(k_i\wedge l_i)\leq k_i$ so this is equivalent to  $(\tau_{k_i}\circ\alpha_{l_i})(b,*) = 1$.
 
 Thus $(\tau\circ\varepsilon_F)_h(b,\alpha) = 1$ if and only if $(\varepsilon_G\circ\Delta_{\inf}\Theta_{\inf}(\tau))_h(b,\alpha) = 1$.

In analagous ways we can show that the required triangles for an adjunction are also satisfied and thus $\varepsilon$ is a natural-isomorphism from which it follows that $\Delta_{\inf}$ is the left adjoint to $\Theta_{\inf}$.\qed

\section{The Equivalences}

The relationships constructed so far are represented in the (not necessarily commutative) diagram below.  It is our goal this section to create three equivalences horizontally as indicated by the dashed arrows.

$$\bfig
\square(0,500)|allb|/>.>>`>->`>->`/<850,500>[\shv`\Relshv`\Preshv`\Relpreshv;\Delta\sub{\Shv}`\iota`\iota`\Delta_{\pre}]
\square/>..>>`>`>`>..>>/<850,500>[\Preshv`\Relpreshv`\shvd`\Relshvd;`\sim`\sim`\Delta\sub{\Shv}]
\morphism(100,600)<600,350>[`;\Delta_{\inf}]
\morphism(200,600)|b|/<-/<600,350>[`;\Theta_{\inf}]
\morphism(910,600)<0,400>[`;\vee]
\put(47.5,87){{\tiny $\dashv$}}
 \efig$$

 We begin this process by showing that the image of $\Theta_{\inf}$ is equivalent to the category of sheaves on $\h$.  This is done in two steps; first it is shown that a presheaf $F$ is a sheaf if and only if $\eta_F$ is an isomorphism,  then we show that $\eta\sub{\Theta(F)}$ is an isomorphism.  This result strongly uses the idea that a singleton morphism $\alpha$ represents the amalgamation of a matching family.  Thus, if $F$ is not a sheaf, the construction $\Theta_{\inf}\Delta_{\inf}(F)$ is adjoining to $F$ all the amalgamations needed to turn $F$ into a sheaf. So $\Theta_{\Shv}\Delta_{\inf}(F)$ will be an associated sheaf functor and in addition $\Delta_{\inf}\Theta_{\pre}$ will be an associated relational-sheaf functor.
 

\begin{lm}\label{equ}\em A presheaf $F$ is a sheaf if and only if $\eta_F$ is an isomorphism.
\end{lm}
\pf Assume that $F$ is a sheaf and let $\alpha:F_h\ra \Delta(F)$ be a singleton.  The family of elements of $\h$ determined by $\alpha$,  $\langle \alpha(x)\rangle_x$,  is a cover of $h$. This follows since $\alpha$ is a symmetric monomorphic map, so $h = \vee(\alpha^{\circ}(x)\wedge \alpha(x))$.  In addition the $x\in X$ are a matching family in $F$ for the cover since by definition whenever $\Delta(F)(k)(x,z) = 1$ then $x_{|_k} = z_{|_k}$.

Thus there is a unique amalgamation $y\in F(h)$ such that  $ x_{|_{\alpha(x)}} = y_{|_{\alpha(x)}}$.  Using this we have

\vsc{-1}

\begin{eqnarray*}
\alpha_k(x,*) = 1 & \Rightarrow &k\leq\vee\{l\ |\ \alpha_l(x,*) = 1\}\\
& \Rightarrow & k \leq \vee\{l\ |\ x_{|_l} = y_{|_l}\}\\
& \Leftrightarrow & \Delta(F)(k)(x,y) = 1\\
&\Leftrightarrow & \alphasingle{y}{k}(x,*) = 1
\end{eqnarray*}

So $\alpha\leq\alphasingle{y}{}$ and since $\Relshv$ is a category it must be the case the $\alpha$ equals $\alphasingle{y}{}$
and thus $\eta_{F,h}$ is a bijection between the sets $F(h)$ and $\Theta\Delta(F)(h)$.

Now assume that $\eta_F$ is an isomorphism and let $\langle h_i\rangle_{i\in I}$ be a cover of $h\in\h$ and $\langle x_i\rangle_{i\in I}$ a matching family for the cover.
Define $\alpha:F_h\Ra F$, by

\begin{center}$\alpha(x) = \left\{\begin{array}{cc} h_i & \hbox{ if } x = x_i\\
							\bot & \hbox{otherwise}
							 \end{array}\right.$\end{center}
							 
It is straight forward to show that $\alpha$ is a singleton on $F$.  Since $\eta$ is an isomorphism there is a $y$ such that $\alpha = \alphasingle{y}{} = \Delta F(-,y)$

Thus $\alpha^y_{h_i}(x_i,*) = 1$ if and only if $\Delta(F)_{h_i}(x_i,y) = 1$ if and only if $x_{i_{|_{h_i}}} = y_{|_{h_i}}$ and hence every matching family of a cover of $h$ has a unique amalgamation $y$ and therefore $F$ is a sheaf. \qed

\begin{thm}\em Let $\h$ be a Heyting algebra, then the category of sheaves of $\h$ is equivalent to the category of relational-sheaves of $\h$.
\end{thm}
\pf The only thing we have left to show is that the image of $\Theta$ is contained within the category of sheaves.  Let $F$ be a relational-sheaf.  We will show that
$\arr{\Theta(F)}{\eta_{\Theta(F)}}{\Theta\Delta\Theta(F)}{}$ is an isomorphism and thus $\Theta(F)$ is a sheaf.  Let $\arr{F_h}{A}{\Delta\Theta(F)}{}$ be a singleton.  We want to find a singleton $\arr{F_h}{\alpha}{F}{}$ such that $A = \Delta\Theta(F)( - ,\alpha)$  (in other words $A = A^{\alpha}$).

Define $\alpha$ by\quad $\ds\alpha_h(x,*) = 1\ \hbox{ if }\ h \leq \bigvee_{\gamma}\Big\{ k\wedge l\ |\ \gamma\sub{k}(x,*)=1\hbox{ and } A\sub{l}(\gamma, *)=1\Big\}$ where  the supremum is taken over all singleton morphisms $\gamma: F_h\ra F$.  To show that $\alpha$ is a singleton morphism we need to show that $F\alpha = \alpha = \alpha F_h$ and that $\alpha$ is a symmetric map.

\begin{eqnarray*}
F\alpha_h(x,*) = 1 & \Leftrightarrow & h\leq\bigvee_y\left\{ k\wedge l \begin{array}{rcl}&|& F_k(x,y) = 1\hbox{ and }\\&|&\\
&|&\ds l\leq \bigvee_{\gamma}\{ r\wedge s\ |\ \gamma_{r}(y,*) = 1\hbox{ and } A_{s}(\gamma,*) = 1\}\end{array}\right\}\\
&&\\
& \Leftrightarrow & h\leq \bigvee_{y,\gamma}\left\{k\wedge r\wedge s\begin{array}{rcl} &|& F_k(x,y) = 1\hbox{ and } \gamma_{r}(y,*) = 1\\ 
&|&\\
&|&\hbox{ and } A_{s}(\gamma,*) = 1\end{array}\right\}\\
&&\\
& \Leftrightarrow & h \leq \bigvee_{\gamma}\left\{k\wedge s\begin{array}{rcl} &|& \ds k\leq \bigvee_y\{k'\wedge r\ |\ F_{k'}(x,y)=1\hbox{ and } \gamma_r(y,*) =1\}\\
&|& \hbox{ and } A_s(\gamma,*)=1 \end{array}\right\}\\
&&\\
& \Leftrightarrow & h\leq \bigvee_{\gamma}\Big\{r\wedge s\ |\ \gamma_r(x,*)=1\hbox{ and }A_s(\gamma,*) = 1\Big\}\\
& \Leftrightarrow & \alpha_h(x,*) = 1
\end{eqnarray*}

So $F\alpha = \alpha$ and in a very similar way $\alpha F_h = \alpha$.  Next to show that $\alpha$ is a map we have

\vsc{-.5}

\begin{eqnarray*}
\alpha\alpha^{\circ}_h(x,y) = 1 & \Leftrightarrow & h\leq \bigvee\left\{k\wedge l\begin{array}{rcl}  & | &\ds k\leq \bigvee_{\gamma}\left\{r\wedge s\begin{array}{rcl} &|&\ds \gamma_r(x,*)=1  \\
&|&\hbox{ and }A_s(\gamma, *) = 1\end{array}\right\}\\
&|&\hbox{ and }\\
& | & \ds l\leq \bigvee_{\xi}\left\{ u\wedge v\begin{array}{rcl} &|& A_v^{\circ}(*, \xi) = 1\\
&|&\hbox{ and } \xi_u^{\circ}(*,y)=1\end{array}\right\} \end{array} \right\}\\
&&\\&&\\
& \Leftrightarrow & h\leq \bigvee_{\gamma,\xi}\left\{\begin{array}{ rcl} & | & \gamma_r(x,*) = 1 \hbox{ and }\\
r\wedge m\wedge v& | & \ds m\leq\{s\wedge u\ |\ A_sA_u^{\circ}(\gamma,\xi) = 1\}\\
& | &  \hbox{ and }\xi^{\circ}_u(*,y) = 1\end{array}\right\}\\
&&\\&&\\
& \Rightarrow & h\leq \bigvee_{\gamma,\xi}\left\{r\wedge m\wedge v\begin{array}{rcl} & | & \gamma_r(x,*)=1\\
& | & \hbox{ and }\Theta\Delta(F)(m)(\gamma,\xi)=1\\
& | & \hbox{ and }\xi^{\circ}_v(*,y) = 1 \end{array}\right\}\\
&&\\&&\\
& \Leftrightarrow & h\leq \bigvee_{\gamma,\xi}\left\{ r\wedge m\wedge v\begin{array}{rcl} &|&\quad\quad \gamma_r(x,*) = 1 \\
&|&\hbox{ and } (\gamma^{\circ}\xi)_m(*,*) = 1\\
&|&\hbox{ and } \xi^{\circ}_v(*,y) = 1
\end{array}\right\}\hbox{ lemma \ref{lm4}}\\
&&\\
& \Rightarrow & (\gamma\gamma^{\circ}\xi\xi^{\circ})_h(x,y) = 1\\
& \Rightarrow & \hspace*{-.3cm}F(h)(x,y) = 1\hbox{ since }\gamma\gamma^{\circ},\xi\xi^{\circ}\leq F\hbox{ for  singleton morphisms}
\end{eqnarray*}
 And thus $\alpha\alpha^{\circ} \leq F$ and similarly $F_h \leq \alpha^{\circ}\alpha$. So $\alpha$ is a symmetric map.    To show that $A = \Delta\Theta(F)(\alpha,-) = A^{\alpha}$ observe that

\begin{eqnarray*}
A_h(\beta,*) = 1 & \Leftrightarrow & \Big(\Theta\Delta(F)A\Big)_h(\beta,*)=1\\
& \Leftrightarrow & h\leq\bigvee_{\gamma}\Big\{k\wedge l\ |\ \Theta\Delta(F)_k(\beta,\gamma) = 1 \hbox{ and } A_l(\gamma,*) = 1 \Big\}\\
& \Leftrightarrow & h\leq \bigvee_{\gamma}\Big\{ k\wedge l\ |\ \beta^{\circ}\gamma_k(*,*) = 1 \hbox{ and } A_l(\gamma,*)=1\Big\} \hbox{ by lemma \ref{lm4} }\\
& \Leftrightarrow & h\leq \bigvee_{x}\left\{ r\wedge s\begin{array}{rcl} &|& \beta^{\circ}_r(*,x) = 1\hbox{ and }\\
& | & \ds s\leq \bigvee_{\gamma}\{k\wedge l\ |\ \gamma_k(x,*) = 1 \hbox{ and } A_l(\gamma,*)=1\}\end{array}\right\} \\
& \Leftrightarrow & h\leq\bigvee_x\Big\{k\wedge l\ |\ \beta^{\circ}_k(*,x) = 1\hbox{ and } \alpha_l(x,*) = 1\Big\}\\
& \Leftrightarrow & \Theta\Delta(F)_h(\beta, \alpha) \hbox{ by lemma \ref{lm4} }\\
& \Leftrightarrow & A^{\alpha}(\beta, *) = 1
\end{eqnarray*}

And our singleton morphism $A$ is the representable morphism $A^{\alpha}$ which tells us that $\eta\sub{\Theta(F)}$ is an isomorphism, hence the category of sheaves for $\h$ is equivalent to the category of relational-sheaves.\qed

\vsc{-.5}

An immediate consequence is that the category of sheaves on the Heyting algebra of down-closed subsets of $\h$ ($\shvd$), is equivalent to the category of relational-sheaves on $\D\h$, ($\Relshvd$).   Combining this with the two comparisons theorems shows that the category of presheaves, $\Preshv$\hspace{-5pt}, is equivalent to the category of relational-sheaves, $\Relpreshv$.  

We can now fill out our diagram of adjuctions and equivalences by setting:

\begin{itemize}

\item {Equivalence of presheaves and relational-presheaves:} 

 $\Delta_{\pre} = \Phi\circ\Theta_{\Shv}\circ\Gamma$ and $\Theta_{\pre} = \Lambda\circ\Delta_{\Shv}\circ\Psi$.

\item { Associated sheaf functor:}  

$a_{\Shv} = \Theta_{\Shv}\circ\Delta_{\inf}:\Preshv\ra\shv$ which is a left adjoint to the functor

$\Theta_{\inf}\circ\Delta_{\Shv}:\shv\ra\Preshv$.

\item { Associated relational-sheaf functor:} 

$a_{\rel} = \Delta_{\inf}\circ\Theta_{\pre}:\map^*\Big(\Relpreshv\Big)\ra\Big(\Relshv\Big)$ which is left adjoint to the functor
$\Delta_{\pre}\circ\Theta_{\inf}:\map^*\Big(\Relshv)\ra\Big(\Relpreshv\Big)$.

\end{itemize}

$$\bfig
\square(0,500)|alla|/->`>->`>->`->/<950,500>[\shv`\map^*(\Relshv)`\Preshv`\map^*(\Relpreshv);\Delta_{\Shv}`\Theta\Delta`\Delta\Theta`\Delta_{\pre}]
\square|allb|/`>`>`->/<950,500>[\Preshv`\map^*(\Relpreshv)`\shvd`\map^*(\Relshvd);`\Gamma`\Psi`\Delta_{\Shv}]
\morphism(150,615)|a|/->/<700,300>[`;\Delta_{\inf}]
\morphism(175,555)|b|/<-/<700,300>[`;\Theta_{\inf}]
\morphism(200,950)|b|/<-/<450,0>[`;\Theta_{\Shv}]
\morphism(300,50)|a|/<-/<270,0>[`;\Theta_{\Shv}]
\morphism(130,450)|b|/<-/<470,0>[`;\Theta_{\pre}]
\morphism(50,600)<0,400>[`;a\sub{\Shv}]
\morphism(1000,600)<0,400>[`;a\sub{\rel}]
\morphism(50,105)<0,400>[`;\Lambda]
\morphism(1000,100)<0,400>[`;\Phi]
\put(50,81.5){\tiny $\bot$}
\put(45,54){\tiny $\sim$}
\put(113,85){{\tiny $\vdash$}}
\put(0.5,87){{\tiny $\vdash$}}
\put(45,113.5){\tiny $\sim$}
\put(48,2){{\tiny $\sim$}}
\put(0.5,30){{\tiny $\sim$}}
\put(111.5,30){{\tiny $\sim$}}
 \efig$$

It follows that every presheaf $F$ on a Heyting algebra is a symmetric idempotent order-preserving lax-semifunctor $F:\h^{co}\ra\rel$ and $F$ is a sheaf if the associated relational-presheaf preserves all infima and is an idempotent using the appropriate construction of the composite.   The associated relational-sheaf functor assigns to a relational-presheaf $F$ the singleton-morphisms that convert the associated down-closed subsets associated to every pair of elements into principle down-closed subsets and thus converting a relational-presheaf into a relational-sheaf.

 In Garraway\cite{wdg1} and Stubbe\cite{s1,s2}, the category of $\q$-valued sets were constructed and in Garraway\cite{wdg1,wdg2} an equivalence between this category of $\q$-valued sets and a generalised notion of sheaves for $\q$ was created.  Then in Garraway\cite{wdg3} an equivalence between $\q$-valued sets and relational-sheaves on $\q$ was also constructed.  For future work, building and extending these results, we anticipate constructing a similar series of equivalences and adjunctions between the categories of sheaves on an involutive quantaloid $\q$ and the category of relational-sheaves on $\q$.  The long term goal is to define the appropriate notion of relational-presheaf  for a site of a category and build the equivalences between the categories.

\vsc{-1}

William Dale Garraway, 

Eastern Washington University,Cheney Washington U.S.A.

  dgarraway@ewu.edu

\end{document}